\documentclass[11pt,dvipsnames,reqno]{amsart}
\usepackage[T1]{fontenc}
\usepackage{amsaddr}
\usepackage{yhmath, mathrsfs, amsthm, amsmath, amssymb, amsfonts, enumerate, lipsum, appendix, mathtools, bm}
\usepackage{mathrsfs}
\usepackage{pst-node}
\usepackage{bbold}
\usepackage[mathscr]{euscript}
\usepackage[normalem]{ulem} 
\usepackage{graphicx}
\usepackage{pgf,tikz}
\usepackage{tkz-euclide}
\usepackage{graphicx}
\usepackage{subcaption}
\usetikzlibrary{shapes.geometric}
\usetikzlibrary{arrows,hobby}
\usepackage{enumitem}
\usepackage[english]{babel}

\allowdisplaybreaks

\usepackage[sort,numbers]{natbib}
\usepackage{hyperref}
\hypersetup{
	colorlinks=true,
	linkcolor=blue,
	filecolor=magenta,      
	urlcolor=red,
	citecolor=red,
}
\usepackage[margin=1in]{geometry}

\usepackage{orcidlink}
\newtheorem{theorem}{Theorem}[section]
\newtheorem{lemma}[theorem]{Lemma}
\newtheorem{proposition}[theorem]{Proposition}
\newtheorem{definition}[theorem]{Definition}
\theoremstyle{definition}

\newtheorem{remark}[theorem]{Remark}
\numberwithin{equation}{section}
\usepackage[hyperpageref]{backref}

\newcommand*\rd{\mathbb{R}^d}

\newcommand{\al} {\alpha}
\newcommand{\pa} {\partial }
\newcommand{\be} {\beta}

\newcommand{\De} {\Delta}

\newcommand{\Om} {\Omega}
\newcommand{\la} {\lambda}
\newcommand{\ga} {\gamma}

\newcommand{\no} {\nonumber}
\newcommand{\noi} {\noindent}

\newcommand{\ra} {\rightarrow}

\def\dx{\,{\rm d}x}
\def\dy{\,{\rm d}y}

\def\ds{\,{\rm d}s}

\def\C{{\mathcal C}}

\def\rd{{\mathbb R}^d}
\def\({{\Big(}}
\def\){{\Big)}}
\def\cc{{\C_c^\infty}}

\DeclareMathOperator{\diam}{diam}

\def\Xint#1{\mathchoice
{\XXint\displaystyle\textstyle{#1}}%
{\XXint\textstyle\scriptstyle{#1}}%
{\XXint\scriptstyle\scriptscriptstyle{#1}}%
{\XXint\scriptscriptstyle\scriptscriptstyle{#1}}%
\!\int}
\def\XXint#1#2#3{{\setbox0=\hbox{$#1{#2#3}{\int}$ }
\vcenter{\hbox{$#2#3$ }}\kern-.6\wd0}}

\def\dashint{\Xint-}
\makeatletter
\newcommand{\vo}{\vec{o}\@ifnextchar{^}{\,}{}}

\DeclarePairedDelimiter\abs{\lvert}{\rvert}%
\DeclarePairedDelimiter\norm{\lVert}{\rVert}%
\def\wps{{W_0^{s,p}(\Om)}}

\title[Gradient regularity estimates]{Lipschitz potential estimates for diffusion with jumps}
\author[N. Biswas and H. Prasad]{Nirjan Biswas$^\dagger$ \and Harsh Prasad}
\address{\rm 
TIFR Centre For Applicable Mathematics Bengaluru\\ Post Bag No 6503, Sharada Nagar, Bengaluru 560065, India}
\email[Nirjan]{nirjan22@tifrbng.res.in} \email[Harsh]{harsh@tifrbng.res.in}
\thanks{$^\dagger$Corresponding author}
\begin{document}
\begin{abstract}
    For $p \in (1, \infty)$ and $s \in (0,1)$, we consider the following mixed local-nonlocal equation 
    \begin{equation}\tag{P}
    \begin{aligned}
          - \Delta_p u +  (-\Delta_p)^s u = f \; \text{in} \; \Omega,
    \end{aligned}
    \end{equation}
    where $\Omega \subset \mathbb{R}^d$ is a bounded domain and the function $f \in L_{loc}^1(\Omega)$. Depending on the dimension $d$, we prove gradient potential estimates of weak solutions to (P) for the entire ranges of $p$ and $s$. 
    As a byproduct, we recover the corresponding estimates in the purely diffusive setup, providing connections between the local and nonlocal aspects of the equation. Our results are new, even for the linear case  $p=2$.
    \end{abstract}
\maketitle
\section{Introduction and Main Result}

The main focus of this paper is to examine Lipschitz potential estimates for a nonlinear jump-diffusion process, which is described by the nonhomogenous equation:
\begin{equation}\label{main-eqn}
\begin{aligned}
- \Delta_p u + (-\Delta_p)^s u = f \; \text{in} \; \Omega,
\end{aligned}
\end{equation}
where $\Delta_p$ is the $p$-Laplace operator defined as $\Delta_p u := \text{div}( \abs{ \partial u}^{p-2} \partial u)$ and $(-\Delta_p)^s$ is the fractional $p$-Laplace operator defined as 
\begin{equation*}
    (-\Delta _p)^s u (x) := \textrm{P.V.}\int _{\mathbb{R}^d} \frac{|u(x)-u(y)|^{p-2}(u(x)-u(y))}{|x-y|^{d+sp}} \,\mathrm{d}y
\end{equation*}
where P.V. stands for the principle value. Here, $\Omega$ represents a bounded domain in $\mathbb{R}^d$ where the dimension $d \ge 2$. The parameters $p \in (1,\infty)$ and $s \in (0,1)$ denote the integrability exponent and the differentiability parameter, respectively. Building upon the concepts and ideas presented in the studies \cite{KM2014MA, KM2014CVPDE, DM2010CVPDE}, this paper establishes criteria that are essentially optimal for the potential function $f$ to obtain Lipschitz regularity for the solution of \eqref{main-eqn}. Furthermore, by employing the strategy outlined in \cite{KK2020}, the paper also achieves improved potential to get Lipschitz regularity in the borderline cases where $p \ge d$, which can be compared to the results in \cite{BY2019}.

In the linear homogeneous case when $p=2$ and $f=0$, the regularity theory for \eqref{main-eqn} was extensively investigated in several papers \cite{BBCK2009, CKS2011, CKS2010, CK2010, F2009}. These studies used probabilistic and analytic techniques to analyze the problem. Conversely, in the nonlinear homogeneous scenario when $p\neq 2$ and $f=0$, the DeGiorgi-Nash-Moser theory for \eqref{main-eqn} was established in \cite{GK2022}. 
The higher H\"{o}lder regularity of the weak solution to \eqref{main-eqn} was investigated by the authors of \cite{GE2023}. They demonstrated that if $f$ lies in the local Lebesgue space $L_{loc}^q(\Omega)$, where $q > d/p$ for $d \ge p$, and $q \ge 1$ for $d < p$, then the weak solution belongs to the local H\"{o}lder space $C_{loc}^{0, \alpha}(\Omega)$ for $0 < \alpha < \min \left\{(pq-d)/q(p-1), sp/(p-1), 1 \right\}$. Notably, in the absence of the forcing term $f$ and under the restrictions $p \ge 2$, and $sp< p-1$, they established that the gradient of the weak solution to \eqref{main-eqn} is locally H\"{o}lder continuous.
In \cite{FM2022}, the authors further improved all these results. They proved that when $f$ belongs to $L^d(\Omega)$, the weak solution to \eqref{main-eqn} exhibits globally almost Lipschitz regularity. They also showed that if $f \in L_{loc}^q(\Omega)$ with $q>d$, then the gradient of the weak solution to \eqref{main-eqn} is locally Hölder continuous.

The borderline Lipschitz regularity result for the Laplace operator was first obtained in \cite{S1981}. In this work the author proved that $u \in L^1(\rd)$ and $\partial u \in L^{d, 1}(\rd; \rd)$ implies $u$ is continuous on $\rd$. This significant finding, combined with the well-known Calderón–Zygmund theory, allows us to conclude that if $\De u \in L^{d, 1}(\rd)$, then $\partial u$ is continuous on $\rd$. Subsequently, in \cite{KM2014CVPDE}, the authors established a nonlinear Stein theorem for weak solutions to a generalized $p$-Laplacian system with coefficients that are not necessarily constant but possess Dini continuity. Later, the nonlinear Stein theorem for the inhomogeneous quasilinear system in differential forms was studied in \cite{S2019CVPDE}.

The journey of pointwise estimates for weak solutions begins with the seminal works of \cite{KP1992, KM1994}. In these papers, utilizing the concept of Wolff potential, the authors established pointwise estimates for solutions to the equation $- \text{div} \, A(x, \partial u) = \mu$ in $\Omega$, where $\mu$ is a nonnegative finite Radon measure  and $A(x, \partial u) \approx \abs{\partial u}^p$.
In \cite{M2011, DM2010}, the authors presented the pointwise gradient estimates of solutions to the general quasilinear equation $-\text{div} \, a(x, \partial u) = \mu$ in $\Omega$, where $\mu$ represents a finite Radon measure, and the vector field $a: \rd \ra \rd$ satisfies specific ellipticity conditions. In the linear case, when $a(z) = z$, this equation simplifies to the Poisson equation. These estimates were obtained using the linear Riesz potentials while assuming that the solution belongs to $C_{loc}^{1, \alpha}(\Omega)$ for a certain $\alpha \in (0,1)$, as stated in \cite{M2011}, and  belongs to $C^1(\Omega)$, as stated in \cite{DM2010}.
In \cite{KM2013}, the authors showed that if $u$ belongs to the Sobolev space $W^{1,p}(\rd)$ and serves as a local weak solution to the equation $- \Delta_p u = \mu$ in $\rd$, with $\mu$ representing a Borel measure with locally finite mass, i.e., $\abs{\mu}(K) < \infty$ for every compact open set $ K \subset \rd$, then $\partial u$ is continuous on $\rd$. Moreover, for every Lebesgue point $X \in \mathbb{R}^d$ of $\partial u$, they derived the following classical representation formula in terms of the Riesz potential:
\begin{align*}
    \abs{\partial u(X)}^{ p-1 } \lesssim_{(d,p)} \int_{\rd} \frac{\rm{d} \abs{\mu} (y) }{\abs{ X - y}^{d-1}}.
\end{align*}
Notably, gradient estimates via linear and nonlinear Riesz potentials studied in \cite{DM2010, KM2014, KM2013, M2011} turn out to be sharper than the Wolff potential bounds in the degenerate case \cite{DM2010, KM2014, FM2021, F2022, BM2020}. 
Several other noteworthy works in this field include \cite{S2019CVPDE, KM2016, KM2014JEMS, BCDKS2018, KM2018, DM2010CVPDE, DM2011, AC, BKH, S}.

In this paper, we obtain the pointwise gradient estimate for the weak solutions to \eqref{main-eqn} in terms of the continuous Dini sum of $f$, as presented in the following theorem.  
\begin{theorem}\label{thm:main}
   Let $f \in L_{loc}^1(\Omega)$, and let $u$ be a weak solution to \eqref{main-eqn}. Suppose $X \in \Omega$ is a Lebesgue point of $\partial u$ and the radius $\rho \in (0,1)$ such that $B_{10\rho}(X) \subset \Omega$. Additionally, when $d=p$, we choose $\rho >0$ small enough so that  $C(d) ( \diam B_{4 \rho} )^d \le 1$, for $C(d)$ as in the John-Nirenberg embedding (Proposition \ref{John-Nirenberg}). Then the following estimate holds: 
    \[
    |\partial u (X)| \leq C \left( \left( \dashint_{B_{\rho}(X)} \left( \abs{\partial  u(x)}^p  + 1 \right) \, \dx \right)^{\frac{1}{p}} +  M(f, B_{4\rho}(X))^{\max\{\frac{1}{p-1},1\}} + T \rho^{ \beta } \right),
    \]
    where $M$ is a potential defined in Definition \ref{defn:pot}, $C= C(d,s,p)$ and $\beta = \beta (d,s,p)$ are positive universal constants, and $T = T(d,s,p, \norm{\partial  u}_{L^p(\Omega)}, \norm{u}_{W^{s,p}(\rd)})$ is a positive constant capturing the jump process.
\end{theorem}

To provide a concise explanation, we choose to focus on a constant coefficient operator. However, it is important to note that our results can be extended to more general structures for both the diffusive and jump terms, provided that appropriate growth conditions are satisfied. The conditions specified in \cite[(1.5)-(1.7) and (1.13)-(1.14)]{FM2022} ensure that the growth behaviour of these operators is compatible with our analysis. 

\begin{remark}
    To recover the results for the diffusion process, i.e. to recover the potential estimates for 
    \begin{equation*}
        \begin{aligned}
            - \Delta_p u = f \; \text{ in } \; \Omega,
        \end{aligned}
    \end{equation*}
    we simply set $T = 0$ in Theorem \ref{thm:main}.
\end{remark}

The remaining sections of the paper are organized as follows. In Section \ref{section-2}, we gather all the preliminary results. Section \ref{section-3} establishes the local comparison estimates between the weak solution to equation \eqref{main-eqn} and its $p$-harmonic replacement, considering the dimension $d$ and the values of $p$. Finally, in Section \ref{section-4}, we present the proof of the main theorem.

\section{Definitions and Auxiliary Results}\label{section-2}
This section collects the preliminary results required in this paper. Initially, we provide a list of notations used throughout the paper.

\noi \textbf{Notation:}
\begin{itemize}
    \item We denote the dimension of the ambient space by $d$.
    \item The ball of radius $r$ centered at $x_0$ is denoted as  $B_r(x_0)$.
    \item We denote the volume of the unit ball in $\rd$ by  $\omega_d$.
    \item For $q \in (1, \infty)$, the conjugate exponent is denoted as $q'= \frac{q}{q-1}$.
    \item We denote the Sobolev exponent $p^* = \frac{dp}{d-p}$ when $d>p$, and $p^*= \infty$ when $d \le p$.
    \item We denote the average integral of a function $h$ over a set $A$: 
    \begin{align*}
         (h)_A := \dashint_A h(x) \, \dx = \frac{1}{|A|} \int_{A} h(x) \, \dx.
    \end{align*}
    \item For $q>0$, we define $|f|^{q}(B_r(x_0)) := \int_{B_r(x_0)} \abs{f(x)}^q \, \dx$.
    \item  For an $L^p$ function $v$ we let $\partial h$ denote its distributional gradient. 
    \item By a universal constant, we mean a constant that depends only on the dimension $d$, the integrability exponent $p$ and the differentiability parameter $s$. 
    \item We use $\lesssim_{(a,b,c)}$ to denote an inequality with a constant depending on $a, b, c$.
\end{itemize}
Now we define the following function spaces. 
\begin{definition}[Sobolev space]
    For $p \in (1, \infty)$ and an open set $\Om \subset \rd$, the Sobolev space is defined as
    \begin{align*}
        W^{1, p}(\Om) := \left\{ u \in L^p(\Omega) : \abs{\partial  u} \in L^p(\Omega) \right\},
    \end{align*}
with the norm $\norm{u}_{W^{1,p}(\Omega)}:= \norm{u}_{L^p(\Omega)} + \norm{\partial  u}_{L^p(\Omega)}$. The subspace $W_0^{1,p}(\Om)$ is the closure of $\C_{c}^{\infty}( \Omega )$ with respect to $\norm{u}_{W^{1,p}(\Omega)}$. 
\end{definition}

\begin{definition}[Fractional Sobolev space]
For $s \in (0,1)$ and $p \in (1, \infty)$, the fractional Sobolev space is defined as
$$ W^{s, p}(\rd):=\left\{u\in L^p(\rd):[u]_{s,p} < \infty \right\},$$
with the fractional Sobolev norm $\norm{u}_{W^{s,p}(\rd)} := \|u\|_{L^p(\rd)} + [u]_{s,p}, $
where
$$ [u]_{s,p}^p :=\iint \limits_{\rd \times \rd}\frac{|u(x)-u(y)|^p}{|x-y|^{d+sp}}\,\dx\dy, $$
is the Gagliardo seminorm.   
\end{definition}

 Now we define the notion of the weak solution to the problem \eqref{main-eqn}.

 \begin{definition}[Weak solution]
     Let $\Om \subset \rd$ be a bounded open set, and let $f \in L_{loc}^1(\Om)$. A function $u \in W^{1,p}(\Om) \cap W^{s, p}(\rd)$ is said to be weak solution to \eqref{main-eqn}, if the following identity holds for every $\phi \in \cc(\Omega)$:
     \[ \int_{\Omega}  \abs{\partial u}^{p-2} \partial u \cdot \partial \phi \, \dx + \iint \limits_{\rd \times \rd}\frac{|u(x)-u(y)|^{p-2} (u(x) -u(y) (\phi(x) - \phi(y))) }{|x-y|^{d+sp}}\,\dx\dy = \int_{\Om} f \phi \, \dx. \]
 \end{definition}
Note that by the density of $\cc(\Omega)$, the above identity holds for every $\phi \in W_0^{1,p}(\Om) \cap W^{s, p}(\rd)$.
We define the Orlicz space (see \cite{AF2003} for details) for the endpoint comparison estimate.
\begin{definition}[Orlicz space]\label{Orlicz}
For a bounded open set $\Om \subset \rd$  we define the Orlicz space $L_{B}(\Om)$ by 
\begin{align*}
    L_{B}(\Om) := \left\{ u \in L_{loc}^1(\Om) : \norm{u}_{B} < \infty \right\},
\end{align*}
where
\begin{align*}
    \norm{u}_{B} := \inf \left\{ s>0 : \int_{ \Om } B \left( \frac{\abs{u(x)}}{s} \right)  \, \dx \le 1 \right\}
\end{align*}
for $B = A$ or $B = \tilde{A}$, where $A$ and $\tilde{A}$ are given by
\begin{align}\label{A and tilde A}
    A(s) = (1+s)\log(1+s) -s \, \text{ and } \, \tilde{A}(s) =  \exp(s) -s -1. 
\end{align}
\end{definition}
The endpoint potential space is obtained as the following Lorentz-Zygmund space. 
\begin{definition}[Lorentz-Zygmund space]
For a bounded open set $\Omega \subset \rd$, the Lorentz-Zygmund space $L \log L(\Omega)$ is defined as: 
\begin{align*}
    L \log L(\Omega) := \left\{ f \in L_{loc}^1(\Omega) : \abs{f}_{L \log L(\Omega)} < \infty \right\},
\end{align*}
where \[\abs{f}_{L \log L(\Omega)} := \int_{\Omega } \abs{ f } \log \left( e + \frac{\abs{f}}{\int_{ \Omega} \abs{f} \, \dx} \right) \, \dx. \]
\end{definition}
In the following definition, we list all the potentials used in this paper. 
\begin{definition}\label{defn:pot}
 For $f \in L^1_{\text{loc}}(\Omega)$ and a ball $B_r(x_0) \subset \Om$ we define the discrete potentials as follows:
 \[F(f, B_r(x_0)) :=\begin{cases*}
    \displaystyle \left(\frac{|f|^{ {(p^*)}'}(B_{r}(x_0))}{{r}^{d-{(p^*)}'}}\right)^{\frac{1}{ {(p^*)}'}},&  if  $d>p$ and $p < 2$;\\
    \displaystyle \left(\frac{|f|^{ {(p^*)}'}( B_{r}(x_0) )}{{r}^{d- {(p^*)}'}}\right)^{\frac{1}{ {(p^*)}'(p-1)}}, & if $p \geq 2 \text{ and } d \neq p$;\\
    \displaystyle \left( {r}^{1-d}|f|_{L\log L(B_{r}(x_0))}\right)^{\frac{1}{d-1}},& if $d=p$.
    \end{cases*}\]
We also define their continuous Dini sum counterpart as follows:
\begin{align*}
     M(f, B_r(x_0)) := \left\{ \begin{array}{ll}
             \displaystyle \int_{0}^r \left( \frac{|f|^{(p^*)'}(B_s(x_0))}{s^{d - (p^*)'}} \right)^{\frac{1}{(p^*)'}} \frac{\ds}{s}, &  \text{ if }  d > p \text{ and } p < 2; \vspace{0.1 cm} \\
             \displaystyle \int_{0}^r \left( \frac{|f|^{(p^*)'}(B_s(x_0))}{s^{d - (p^*)'}} \right)^{\frac{1}{(p^*)'(p-1)}} \frac{\ds}{s}, &  \text{ if }  p \ge 2 \text{ and } d \neq p; \vspace{0.1 cm} \\
             \displaystyle \int_{0}^r \left( s^{1-d} \abs{f}_{ L \log L (B_s(x_0)) } \right)^{\frac{1}{d-1}} \frac{\ds}{s}, &  \text{ if }  d = p.
            \end{array}\right.
 \end{align*}
\end{definition}
In the next definition, we talk about the excess functional. 
\begin{definition}[Excess functional]
    For a measurable set $A \subset \rd$ and $h \in L^1(A)$, the excess functional is defined as
    \begin{align*}
        \mathfrak{E}_p(h, A) := \left( \dashint_{A} \abs{h(x) - (h)_A}^p  \, \dx \right)^{\frac{1}{p}}.
    \end{align*}
\end{definition}

Now, we state all the auxiliary results. We begin with the following algebraic inequalities. For proof, see \cite[Chapter 8]{G2003}.
\begin{lemma}
Let $p \in (1, \infty)$. Then the following hold: 

 \begin{enumerate}[label=\rm (\roman*)]
    \item For any two vectors $A_1, A_2 \in \rd$, the following holds:
    \begin{equation}\label{ineq1}
        \begin{split}
            \left| \abs{A_1}^{\frac{p-2}{2}} A_1 - \abs{A_2}^{\frac{p-2}{2}} A_2 \right|^2 & \lesssim_{(d, p)} \left( \abs{A_1}^2 + \abs{A_2}^2 \right)^{\frac{p-2}{2}} \abs{A_1 - A_2}^2 \\
    & \lesssim_{(d, p)} \left< \abs{A_1}^{p-2}A_1 - \abs{A_2}^{p-2}A_2 , A_1 - A_2 \right>,
        \end{split}
    \end{equation}
where $\left< \cdot \right>$ denotes the product of vectors. 
\item  Let  $A_1, A_2 \in \rd$ be any two vectors in $\rd$. Then 
the following hold
\begin{align}\label{ineq2}
\text{ for } p \ge 2: \abs{A_1 - A_2}^p \lesssim_{(d, p)} \left( \abs{A_1}^2 + \abs{A_2}^2 \right)^{\frac{p-2}{2}} \abs{A_1 - A_2}^2, 
\end{align}
and 
\begin{equation}\label{ineq2.1}
    \begin{split}
        \text{ for } p < 2: \abs{A_1 - A_2}^p & \lesssim_{(d, p)} \left( \abs{A_1}^2 + \abs{A_2}^2 \right)^{\frac{p-2}{2}} \abs{A_1 - A_2}^2  \\
 & \quad + \left( \left( \abs{A_1}^2 + \abs{A_2}^2 \right)^{\frac{p-2}{2}} \abs{A_1 - A_2}^2 \right)^{{\frac{p}{2}}} \abs{A_1}^{\frac{p(2-p)}{2}}.
    \end{split}
\end{equation}
\end{enumerate}
 \end{lemma}

For a function $v \in \wps$, we denote the duality action $((-\De_p)^s u)(v)$ which is defined as 
\begin{align*}
    ((-\De_p)^s u)(v) := \iint\limits_{\rd \times \rd} \frac{\abs{u(x) - u(y)}^{p-2}(u(x) -u(y))(v(x) -v(y))}{\abs{x-y}^{d+sp}} \, \dx \dy. 
\end{align*}
Now we state the following monotonicity property of the fractional $p$-Laplacian (see \cite[Lemma 2.3 and Lemma 2.4]{IaMoSq1}).

\begin{proposition}\label{monotonicity}
     The following estimate holds for every $u, v \in W^{s,p}(\Omega)$:
    \begin{align*}
      & \text{ for } p \ge 2:  [u-v]_{s,p}^p \lesssim_{(p)} ((-\De_p)^s u - (-\De_p)^s v)(u-v). \\
      & \text{ for } p < 2: \frac{[u-v]_{s,p}^2}{\left( [u]_{s,p}^p + [v]_{s,p}^p  \right)^{2-p}} \lesssim_{(p)} ((-\De_p)^s u - (-\De_p)^s v)(u-v).
    \end{align*}
\end{proposition}
 
The next proposition states the Sobolev-Poincar\'{e} embedding (see \cite[Lemma 1.64]{MZ1997}) results on the ball $B_r(x_0)$.

\begin{proposition}\label{sobolev-embed}
    Let $d>p$ and $q \in [1, p^*]$. Then the following hold:
    \begin{align*}
        & \left( \dashint_{B_r(x_0)} \left| \frac{u(x)}{r} \right|^q \, \dx \right)^{\frac{1}{q}} \lesssim_{(d, p)} \left( \dashint_{B_r(x_0)} \abs{\partial  u(x)}^p \, \dx  \right)^{\frac{1}{p}}, \quad \forall \, u \in W_0^{1,p}(B_r(x_0)), \text{ and } \\
        & \left( \dashint_{B_r(x_0)} \left| \frac{u(x) - (u)_{B_r(x_0)}}{r} \right|^q \, \dx \right)^{\frac{1}{q}} \lesssim_{(d, p)} \left( \dashint_{B_r(x_0)} \abs{\partial  u(x)}^p \, \dx  \right)^{\frac{1}{p}}, \quad \forall \, u \in W^{1,p}(B_r(x_0)).
    \end{align*}
\end{proposition}
We state the following John-Nirenberg embedding for BMO spaces (see \cite[Theorem 1.66]{MZ1997}). 
\begin{proposition}\label{John-Nirenberg}
Let $u \in W^{1,1}(B_r(x_0))$. Suppose there exists a constant $M >0$ such that 
\begin{align*}
    \int_{B_r(x_0) \cap B_{\rho}(x_0)} \abs{\partial  u(x)} \, \dx \le M \rho^{d-1}, \text{ for every } B_{\rho}(x_0) \subset \rd.
\end{align*}
Then there exist positive constants $\sigma_0$ and $C$ depending only on $d$ such that the following inequality holds whenever $\sigma < \sigma_0 | B_r(x_0) | (\diam B_r(x_0))^{-d}$:
\begin{align*}
    \int_{ B_r(x_0) } \exp \left( \frac{\sigma}{M} \abs{ u - (u)_{B_r(x_0)}} \right) \le C ( \diam B_r(x_0) )^d.
\end{align*}
\end{proposition}

Next, we state the Morrey embedding (see \cite[Theorem 1.62]{MZ1997}).

\begin{proposition}\label{Morrey}
    Let $d < p$ and  $u \in W_0^{ 1, p}(\Om)$. Then $u \in C^{0, 1 - \frac{d}{p}}(\overline{ \Omega })$ and the following holds: 
    \begin{align*}
        \sup_{ x \in \Omega  } \abs {u (x) } \lesssim_{(d,p)} \abs{ \Omega }^{\frac{1}{d}} \left(  \dashint_{ \Omega } \abs{ \partial  u(x) }^p \, \dx \right)^{\frac{1}{p}}.
    \end{align*}
\end{proposition}
We need the following equivalence of norms in the Orlicz space. 

 \begin{proposition}\label{equivalence}
    Let $\Omega \subset \rd$ be a bounded open set and $A$ be as given in Definition \ref{Orlicz}. Then there exist constants $k_1, k_2>0$ such that following holds:
    \begin{align*}
        k_1 \norm{f}_{L_A(\Omega)} \le \abs{f}_{L \log L(\Omega)} \le k_2 \norm{f}_{L_A(\Omega)}.
    \end{align*}
 \end{proposition}
 \begin{proof}
     Consider the function $J(s) = s\log(1+s)$ for $s>0$. Then we have 
     \begin{align*}
         \lim_{s \ra \infty} \frac{A(s)}{J(s)} = 1.
     \end{align*}
     Hence using the fact that $|\Omega|$ is finite, we apply \cite[Section 8.4 and Theorem 8.12-(b)]{AF2003}, to get that the norms $\norm{\cdot}_{L_A(\Omega)}$ and $\norm{\cdot}_{L_J(\Omega)}$ are equivalent. Further, by \cite[Lemma 8.6]{IV1999}, $\norm{\cdot}_{L_J(\Omega)}$ and $\abs{\cdot}_{L \log L(\Omega)}$ are also equivalent. Thus the claim follows.  
 \end{proof}
 Next, we state the classical estimates for a constant coefficient homogeneous system. The following proposition is proved in \cite[Theorem 3]{KM2014CVPDE}.
 \begin{proposition}\label{prop:decay1}
    Let $v \in W^{1,p}(B_r(x_0))$ be a weak solution to
    \[
    -\Delta_p v = 0 \,\text{ in } \, B_r(x_0). 
    \]
    Suppose that for some $\la \geq 1$ we have 
    \[
    \dashint_{B_r(x_0)}|\partial  v|^p \lesssim_{(d,p)} \lambda^p.
    \]
    Then there exists an $\al = \al(d,p) \in (0,1)$ such that for any $0<t<1$ we have 
    \[
        \sup_{\frac{1}{2}B_r(x_0)} \, |\partial v| \lesssim_{(d,p)} \lambda \; \text{ and } \; \underset{\frac{t}{2}B_r(x_0)}{\text{osc}} \, \partial v \lesssim_{(d,p)} t^{\alpha}\lambda.
    \]
\end{proposition}

The subsequent decay estimate is proved in \cite[Theorem 3.1]{KM2012}.

\begin{proposition}\label{prop:decay2}
  Let $v \in W^{1,p}(B_r(x_0))$ be a weak solution to
    \[
    -\Delta_p v = 0 \; \text{ in } \; B_r(x_0). 
    \]
    Suppose $\la \geq 1$, $C_0\geq 1$ and $C_1 \geq 1$ are constants. Then for any prescribed decay $\gamma \in (0,1)$ we can find a scale $\sigma' = \sigma'(d,p,C_0,C_1,\gamma) \in (0, \frac{1}{2})$ such that if $\sigma \in (0,\sigma']$ and 
    \begin{align}\label{eq:non-deg}
        \frac{\la}{C_0} \leq \sup_{B_{\sigma r}(x_0)}|\partial v| \leq \sup_{B_{\frac{r}{2}}(x_0)}|\partial v| \leq C_1\la,
    \end{align}
    then
    \begin{align}\label{eq:exdec}
         \mathfrak{E}_p(\partial v,B_{\sigma r}(x_0)) \leq \gamma \mathfrak{E}_p(\partial v, B_{r}(x_0)).
    \end{align}
    
\end{proposition}
The forthcoming lemma establishes a lower bound for $M(f, B_r(x_0))$ in terms of the Dini sum.

\begin{lemma}\label{Dini sum expression}
    Suppose $M(f, B_r(x_0))$ is finite. Let $ \sigma \in (0,1)$ and $r_j = \sigma^j r$. Then the following hold:
   \begin{enumerate}[label=\rm (\roman*)]
       \item Let $d>p$. Then
       \begin{align*}
        & \text{ for } p \ge 2:  \sum_{1 \le j < \infty} \left( \frac{\abs{f}^{(p^*)'}(B_{r_j}(x_0))}{r_j^{d-(p^*)'}} \right)^{\frac{1}{(p^*)'(p-1)}} \lesssim_{(d,p, \sigma)} M(f, B_r(x_0)). \\
        & \text{ for } p < 2:  \sum_{1 \le j < \infty} \left( \frac{\abs{f}^{(p^*)'}(B_{r_j}(x_0))}{r_j^{d-(p^*)'}} \right)^{\frac{1}{(p^*)'}} \lesssim_{(d, \sigma)} M(f, B_r(x_0)).
       \end{align*}
       \item Let $ d =p$. Then 
       \begin{align*}
           \sum_{1 \le j < \infty} \left( r_j^{-(d-1)} \abs{f}_{ L \log L (B_{r_j}(x_0)) }   \right)^{\frac{1}{d-1}} \lesssim_{(\sigma)} M(f, B_r(x_0)). 
       \end{align*}
       \item Let $d<p$. Then 
        \begin{align*}
            \sum_{1 \le j < \infty} \left( \frac{\abs{f}(B_{r_j}(x_0))}{r_j^{d-1}} \right)^{\frac{1}{p-1}} \lesssim_{(d,p, \sigma)} M(f, B_r(x_0)). 
        \end{align*}
   \end{enumerate}
\end{lemma}

\begin{definition}
Let $r>0$ be such that $B_r(x_0) \subset \Om$. Given constants $K_1,K_2 \ge 1$, we say that the pair $(f,\la)$ satisfies the \textbf{finiteness condition}, if the following inequality holds:
\begin{equation}\tag{\textbf{$\Phi$}}\label{H}
    \begin{split}
        K_1^p \left( \dashint_{B_r(x_0)} \left( \abs{\partial  u(x)}^p  + 1 \right) \, \dx \right) + K_2^{p} M(f, B_{r}(x_0))^{\max\{\frac{p}{p-1},p\}}  \le \la^p.
    \end{split}
\end{equation}

\end{definition}

 We now show that if $(f, \la)$ satisfies the finiteness condition, then for every $r>0$, the quantity $(r|f|)^{(p^*)'}(B_r(x_0))$ is estimated by a constant which  depends on $\la, K_2,d,$ and $p$.
\begin{lemma}\label{Dini sum}
   Let $d>p$, $\sigma \in (0,1)$, and $r_j = \sigma^j r$ for some $\sigma \in (0,1)$. Let $d>p$. If $(f, \la)$ satisfies \eqref{H}, then the following holds:
 \begin{align*}
       \left( \dashint_{B_r(x_0)} \abs{r f(x)}^{(p^*)'} \, \dx  \right)^{\frac{1}{(p^*)'}} \lesssim_{(d, p)} \left( \frac{\la}{K_2} \right)^{p-1}.
   \end{align*}
  \end{lemma}
\begin{proof}
For $r>0$ there exists $j \in \mathbb{N}$ and $\rho>0$ such that $2^{-(j+1)}\rho \le r \le 2^{-j} \rho$. Hence using (i) of Lemma \ref{Dini sum expression} we write 
    \begin{align*}
        \left( \dashint_{B_r(x_0)} \abs{r f(x)}^{(p^*)'} \, \dx   \right)^{\frac{ 1 }{(p^*)'}} & \lesssim_{(d)}  \left( \dashint_{B_{2^{-j} \rho}(x_0)} (2^{-j} \rho)^{(p^*)'} \abs{f(x)}^{(p^*)'} \, \dx   \right)^{\frac{ 1 }{(p^*)'}} \\
        & \lesssim_{(d)} \sum_{1 \le j < \infty} \left( \frac{\abs{f}^{(p^*)'}(B_{2^{-j} \rho}(x_0)}{(2^{-j} \rho)^{d-(p^*)'}} \right)^{\frac{ 1 }{(p^*)'}}  \\
        & \lesssim_{(d, p)} M(f, B_r(x_0)) \lesssim_{(d, p)} \left( \frac{\la}{K_2} \right)^{p-1}, 
    \end{align*}
where the last inequality follows using \eqref{H}.
\end{proof}

\section{Comparison estimate}\label{section-3}

In this section, we prove comparison estimates of the weak solution of \eqref{main-eqn} and its $p$-harmonic replacement, in terms of the potential $f$. Throughout this section, we assume $x_0 \in \Omega$ and choose a fixed value of $0 < r < 1$ such that $B_{10 r}(x_0) \subset \Omega$. First, we consider the case $d > p$.

\begin{proposition}\label{comparison estimates d>p}
Let $d > p$. Suppose $u$ is a weak solution to \eqref{main-eqn}. Then there exists $h$ satisfying $- \De_p h = 0 \mbox{ in } B_r(x_0)$, such that the following estimates hold:
\begin{align}\label{h est-1}
    & \text{ for } p \ge 2:  \dashint_{B_r(x_0)} \abs{ \partial (u - h)(x) }^p \, \dx \le C \left( \dashint_{B_{4r}(x_0)} \abs{r f(x)}^{(p^*)'} \, \dx  \right)^{\frac{p'}{(p^*)'}} + T r^{p \be}.
 \end{align}
Further, if $(f, \la)$ satisfies \eqref{H}, then 
\begin{equation}\label{h est-2}
    \begin{split}
        \text{ for } p < 2:  \dashint_{B_r(x_0)} \abs{ \partial (u - h)(x) }^p \, \dx & \le C \left( \dashint_{B_{4r}(x_0)} \abs{r f(x)}^{(p^*)'} \, \dx \right)^{\frac{p}{(p^*)'}} \\
        & \left( \left( \frac{\la}{K_1} \right)^{p(2-p)} + \left( \frac{\la}{K_2} \right)^{p(2-p)}\right) + T r^{p \be},
    \end{split}
\end{equation}
where $C= C(d, p) >0$ and $T = T(d,s,p, \norm{\partial  u}_{L^p(\Omega)}, \norm{u}_{W^{s,p}(\rd)}) >0$ are constants, and the exponent $\be = \be(d,s,p) \in (0,1)$. 
\end{proposition}

\begin{proof}
We divide our proof into two steps. 

\noi \textbf{Step 1:} For $s \in (0,1)$, and $p \in (1, \infty)$, we consider the following mixed local and nonlocal equation: 
\begin{equation}\label{v-eqn}
    \begin{aligned}
        - \De_p v + (-\De_p)^s v &= 0 \mbox{ in } B_{4r}(x_0),\\
        v&=u \mbox{ in } \mathbb{R}^d\setminus B_{4r}(x_0).
    \end{aligned}
\end{equation}
Using the direct variational method, there exists $v \in W^{1, p}(B_{4r}(x_0)) \cap W^{s,p}(\rd)$ satisfying \eqref{v-eqn} weakly.  For brevity, we denote $B_{4r}(x_0)$ as $B_{4r}$. Taking $\phi := u-v \in W_0^{1,p}(B_{4r}(x_0)) \cap W^{s,p}(\rd)$ as a test function in the weak formulations of \eqref{main-eqn} and \eqref{v-eqn}, we have
\begin{equation}\label{est-1}
    \begin{split}
        \dashint_{B_{4r}} \abs{\partial  u}^{p-2} \partial  u \cdot \partial  \phi \, \dx + \frac{1}{\abs{B_{4r}}} ((-\De_p)^s u)(\phi) = \dashint_{B_{4r}} f \phi \, \dx,
    \end{split}
\end{equation}
and 
\begin{equation}\label{est-2}
    \begin{split}
        \dashint_{B_{4r}} \abs{\partial  v}^{p-2} \partial  v \cdot \partial  \phi \, \dx + \frac{1}{\abs{B_{4r}}} ((-\De_p)^s v)(\phi) = 0.
    \end{split}
\end{equation}
In this step, we estimate $\abs{\partial  \phi}^p$ over $B_{4r}$. We separate our proof into two cases.

\noi \textbf{For $p \ge 2$:} In this case, using \eqref{ineq1}, \eqref{ineq2}, and monotonicity property of the fractional $p$-Laplacian (Proposition \ref{monotonicity}) we arrive at
\begin{equation}\label{est-0}
    \begin{split}
         \text{ L.H.S. of } \eqref{est-1} -  \text{ L.H.S. of } \eqref{est-2} & \ge C(d,p) \left( \dashint_{B_{4r}} \abs{\partial  \phi }^p \, \dx  + \frac{1}{|B_{4r}|} [ \phi ]_{s,p}^p  \right) \\
         & \ge C(d,p) \dashint_{B_{4r}} \abs{\partial  \phi }^p \, \dx.
    \end{split}
\end{equation}
Now, using H\"{o}lder's inequality with the conjugate pair $(p^*, (p^*)')$, the Sobolev inequality (Proposition \ref{sobolev-embed}), and then using Young's inequality with the conjugate pair $(p, p')$ we obtain 
\begin{align}
    \left| \dashint_{B_{4r}} f \phi \, \dx \right| & \le \left( \dashint_{B_{4r}} \abs{r f(x)}^{(p^*)'} \, \dx \right)^{\frac{1}{(p^*)'}} \left( \dashint_{B_{4r}}  \left( \frac{\abs{\phi}}{r} \right)^{p^*} \, \dx \right)^{\frac{1}{p^*}} \no \\
    & \label{est-of-f-1} \lesssim_{(d, p)} \left( \dashint_{B_{4r}} \abs{r f(x)}^{(p^*)'} \, \dx \right)^{\frac{1}{(p^*)'}}  \left( \dashint_{B_{4r}} \abs{ \partial  \phi}^p \, \dx \right)^{\frac{1}{p}} \\
    & \label{est-of-f} \le \ga \dashint_{B_{4r}} \abs{ \partial  \phi }^p \, \dx + C( d, p , \ga) \left( \dashint_{B_{4r}} \abs{r f(x)}^{(p^*)'} \, \dx \right)^{\frac{p'}{(p^*)'}}, 
\end{align}
for some $\ga \in (0,1)$. Take $\gamma < \frac{1}{2}$. Therefore, in view of \eqref{est-0}, subtracting \eqref{est-2} from \eqref{est-1}  we get 
\begin{align}\label{ineq 1}
    \dashint_{B_{4r}(x_0)} \abs{ \partial  \phi }^p \, \dx \lesssim_{(d, p)} \left( \dashint_{B_{4r}(x_0)} \abs{r f(x)}^{(p^*)'} \, \dx  \right)^{\frac{p'}{(p^*)'}}.
\end{align}
\noi \textbf{For $p<2$:} In this case, using \eqref{ineq2.1} we write 
\begin{equation}\label{est-3}
    \begin{split}
        \dashint_{B_{4r}} \abs{\partial  \phi }^p \, \dx & \lesssim_{(d,p)} \dashint_{B_{4r}} \left( \abs{\partial  u}^2 + \abs{ \partial  v}^2  \right)^{\frac{p-2}{2}} \abs{\partial  \phi }^2 \, \dx \\
        & \quad +  \dashint_{B_{4r}} \left( \left( \abs{\partial  u}^2 + \abs{\partial  v}^2 \right)^{\frac{p-2}{2}} \abs{\partial  \phi }^2 \right)^{{\frac{p}{2}}} \abs{\partial  u}^{\frac{p(2-p)}{2}} \, \dx.
    \end{split}
\end{equation}
Applying H\"{o}lder's inequality with the conjugate pair $\left( \frac{2}{p}, \frac{2}{2-p} \right)$, we get 
\begin{equation*}
    \begin{split}
        & \dashint_{B_{4r}} \left( \left( \abs{\partial  u}^2 + \abs{\partial  v}^2 \right)^{\frac{p-2}{2}} \abs{\partial  \phi }^2 \right)^{{\frac{p}{2}}} \abs{\partial  u}^{\frac{p(2-p)}{2}} \, \dx \\
        & \le \left( \dashint_{B_{4r}} \left( \abs{\partial  u}^2 + \abs{\partial  v}^2 \right)^{\frac{p-2}{2}} \abs{\partial  \phi }^2 \, \dx \right)^{\frac{p}{2}}
        \left(  \dashint_{B_{4r}} \abs{\partial  u}^p \, \dx \right)^{\frac{2-p}{2}}. 
    \end{split}
\end{equation*}
Since $(f, \lambda)$ satisfies \eqref{H}, we obtain from \eqref{est-3} that  
\begin{equation}\label{est-4}
    \begin{split}
         \dashint_{B_{4r}} \abs{\partial  \phi }^p \, \dx & \lesssim_{(d,p)} \dashint_{B_{4r}} \left( \abs{\partial  u}^2 + \abs{ \partial  v}^2  \right)^{\frac{p-2}{2}} \abs{\partial  \phi }^2 \, \dx \\
         & \quad +  \left( \dashint_{B_{4r}} \left( \abs{\partial  u}^2 + \abs{\partial  v}^2 \right)^{\frac{p-2}{2}} \abs{\partial  \phi}^2 \, \dx \right)^{\frac{p}{2}} \left( \frac{\la}{K_1} \right)^{\frac{p(2-p)}{2}}.
    \end{split}
\end{equation}
Further, using monotonicity property of the fractional $p$-Laplacian in the singular case $1<p<2$ (Proposition \ref{monotonicity}) we have  
\begin{align}\label{est-5}
    0 \le \frac{[\phi]_{s,p}^2}{\left( [u]_{s,p}^p + [v]_{s,p}^p  \right)^{2-p}} \lesssim_{(p)} ((-\De_p)^s u - (-\De_p)^s v)(\phi).
\end{align}
Now we subtract \eqref{est-2} from \eqref{est-1}, and use the estimates \eqref{ineq1}, \eqref{est-of-f}, and \eqref{est-5} to get 
\begin{equation}\label{est-6}
    \begin{split}
        \dashint_{B_{4r}} \left( \abs{\partial  u}^2 + \abs{ \partial  v}^2  \right)^{\frac{p-2}{2}} \abs{\partial  \phi }^2 \, \dx & \le \gamma \dashint_{B_{4r}} \abs{ \partial \phi}^p \, \dx 
        \\
        &\quad + C( d, p , \gamma) \left( \dashint_{B_{4r}} \abs{r f(x)}^{(p^*)'} \, \dx \right)^{\frac{p'}{(p^*)'}}.
    \end{split}
\end{equation}
Using \eqref{ineq1}, \eqref{est-of-f-1}, and \eqref{est-5}, we further estimate the last integral of \eqref{est-4} as: 
\begin{equation}\label{est-7}
    \begin{split}
       & \left( \dashint_{B_{4r}} \left( \abs{\partial  u}^2 + \abs{\partial  v}^2 \right)^{\frac{p-2}{2}} \abs{\partial  \phi}^2 \, \dx \right)^{\frac{p}{2}} \left( \frac{\la}{K_1} \right)^{\frac{p(2-p)}{2}} \\
       & \lesssim_{(d, p)} \left( \dashint_{B_{4r}} \abs{r f(x)}^{(p^*)'} \, \dx \right)^{\frac{p}{2(p^*)'}}  \left( \dashint_{B_{4r}} \abs{ \partial  \phi}^p \, \dx \right)^{\frac{1}{2}} \left( \frac{\la}{K_1} \right)^{\frac{p(2-p)}{2}} \\
       & \le \tilde{\ga} \dashint_{B_{4r}} \abs{ \partial  \phi}^p \, \dx +C(d, p, \Tilde{\ga}) \left( \dashint_{B_{4r}} \abs{r f(x)}^{(p^*)'} \, \dx \right)^{\frac{p}{(p^*)'}} \left( \frac{\la}{K_1} \right)^{p(2-p)},
    \end{split}
\end{equation}
where the last inequality follows using Young's inequality with the conjugate pair $\left ( \frac{1}{2}, \frac{1}{2} \right)$, and $\Tilde{\ga} \in (0,1)$.
Choosing $\ga, \Tilde{\ga} < \frac{1}{4}$ we obtain from \eqref{est-4}, \eqref{est-6}, and \eqref{est-7} that
\begin{equation}\label{est-8}
    \begin{split}
        \dashint_{B_{4r}} \abs{\partial  \phi}^p \, \dx & \lesssim_{(d,p)} \left( \dashint_{B_{4r}} \abs{r f(x)}^{(p^*)'} \, \dx \right)^{\frac{p'}{(p^*)'}}  \\
        & + \left( \dashint_{B_{4r}} \abs{r f(x)}^{(p^*)'} \, \dx \right)^{\frac{p}{(p^*)'}} \left( \frac{\la}{K_1} \right)^{p(2-p)}.
    \end{split}
\end{equation}
Next, we estimate the second integral of \eqref{est-8}. Using the splitting $\frac{1}{p-1} = 1 + \frac{2-p}{p-1}$, we write 
\begin{equation*}
    \begin{split}
        \left( \dashint_{B_{4r}} \abs{r f(x)}^{(p^*)'} \, \dx \right)^{\frac{p'}{(p^*)'}} = \left( \dashint_{B_{4r}} \abs{r f(x)}^{(p^*)'} \, \dx \right)^{\frac{p}{(p^*)'}} \left( \dashint_{B_{4r}} \abs{r f(x)}^{(p^*)'} \, \dx \right)^{\frac{p(2-p)}{(p^*)'(p-1)}}.
    \end{split}
\end{equation*}
Moreover, since $(f, \la)$ satisfies \eqref{H}, using Lemma \ref{Dini sum} we get 
\begin{align*}
    \left( \dashint_{B_{4r}} \abs{r f(x)}^{(p^*)'} \, \dx \right)^{\frac{p(2-p)}{(p^*)'(p-1)}} \lesssim_{(d, p)} \left( \frac{\la}{K_2} \right)^{p(2-p)}.
\end{align*}
Therefore, \eqref{est-8} yields
\begin{align}\label{ineq 2}
    \dashint_{B_{4r}} \abs{\partial  \phi}^p \, \dx \lesssim_{(d, p)} \left( \dashint_{B_{4r}} \abs{r f(x)}^{(p^*)'} \, \dx \right)^{\frac{p}{(p^*)'}} \left( \left( \frac{\la}{K_1} \right)^{p(2-p)} + \left( \frac{\la}{K_2} \right)^{p(2-p)}\right).
\end{align}

\noi \textbf{Step 2:} We consider $h \in W^{1,p}(B_{r}(x_0))$ as a weak solution to the following problem: 
\begin{equation}\label{h-eqn}
    \begin{aligned}
        - \De_p h &= 0 \mbox{ in } B_{r}(x_0),\\
        h&=v \mbox{ on }  \pa B_{r}(x_0).
    \end{aligned}
\end{equation}
Clearly, $v-h \in W_0^{1,p}(B_r(x_0)$. We apply \cite[Lemma 6.2]{FM2022} to obtain 
\begin{align}\label{ineq 3}
    \dashint_{B_r}  \abs{\partial  (v-h)}^p \, \dx \le T r^{p \beta},
\end{align}
where $\beta = \beta(d,s,p) \in (0,1)$, and $T=T(d, p,s,\norm{v}_{L^p(\rd)})$ where the dependence on the norm of $v$ is monotone non-decreasing. Further, we use the energy estimate for \eqref{v-eqn} to get 
\begin{align*}
    \int_{\rd} \abs{v}^p = \int_{B_{4r}} \abs{v}^p + \int_{B_{4r}^c} \abs{u}^p \lesssim_{(d,s)} \left( \norm{u}^p_{W^{s,p}(\rd)} + \int_{B_{4r}} \abs{\partial  u}^p \right). 
\end{align*}
Therefore, \eqref{ineq 3} holds with $T = T(d,s,p, \norm{\partial  u}_{L^p(\Omega)}, \norm{u}_{W^{s,p}(\rd)})$. We split
\begin{align*}
    \dashint_{B_r} \abs{\partial (u - h)}^p \, \dx \lesssim_{(p)} \left( \dashint_{B_r} \abs{\partial (u - v)}^p \, \dx + \dashint_{B_r} \abs{\partial (v - h)}^p \, \dx \right).
\end{align*}
In view of the above inequality, we use \eqref{ineq 1}, \eqref{ineq 2}, and \eqref{ineq 3} to get the required estimates. This completes the proof. 
\end{proof}

In the following proposition, we consider the borderline case $d=p$.

\begin{proposition}\label{comparison estimates for d=p}
Let $d = p$, and $u$ be a weak solution to \eqref{main-eqn}. For $C(d)$ as given in Proposition \ref{John-Nirenberg}, we choose $r>0$ small enough so that  $C(d) ( \diam B_{4r} )^d \le 1$. Then there exists $h$ satisfying $-\De_p h = 0 \mbox{ in } B_r(x_0)$, such that the following estimate holds:
\begin{align}\label{h est-3}
  \dashint_{B_r(x_0)} \abs{ \partial (u - h)(x) }^d \, \dx \le C \left( \dashint_{B_{4r}(x_0)}  r \abs{ f } \log \left( e + \frac{\abs{f}}{\int_{B_r(x_0)} \abs{f} \, \dx} \right) \, \dx  \right)^{\frac{d}{ d-1 }} + T r^{d \be} ,
 \end{align}
where $C = C(d) >0 $ and $T = T(d,s, \norm{\partial  u}_{L^d(\Omega)}, \norm{u}_{W^{s,d}(\rd)}) > 0$ are constants, and $ \be = \be(d,s) \in (0,1)$. 
\end{proposition}

\begin{proof}
    As in the previous proposition, for $s \in (0,1)$, we consider the following mixed local and nonlocal equation: 
\begin{equation}\label{v-eqn critical}
    \begin{aligned}
        - \De_d v + (-\De_d)^s v &= 0 \mbox{ in } B_{4r},\\
        v&=u \mbox{ in } \mathbb{R}^d\setminus B_{4r}.
    \end{aligned}
\end{equation}
For a solution $v \in W^{1, d}(B_{4r}(x_0)) \cap W^{s, d}(\rd)$ of \eqref{v-eqn critical}, we take $\phi := u-v \in W_0^{1,d}(B_{4r}(x_0)) \cap W^{s , d}(\rd)$ as a test function in the weak formulations of \eqref{main-eqn} and \eqref{v-eqn critical} to get
\begin{equation}\label{est-1 critical}
    \begin{split}
        & \dashint_{B_{4r}} \abs{\partial  u}^{d-2} \partial  u \cdot \partial  \phi \, \dx + \frac{1}{\abs{B_{4r}}} ((-\De_d)^s u)(\phi)
         = \dashint_{B_{4r}} f \phi \, \dx, \\
       &  \dashint_{B_{4r}} \abs{\partial  v}^{d-2} \partial  v \cdot \partial  \phi \, \dx + \frac{1}{\abs{B_{4r}}} ((-\De_d)^s v)(\phi)
         = 0. 
    \end{split}
\end{equation}
Therefore, using \eqref{ineq1}, \eqref{ineq2}, and Proposition \ref{monotonicity} we similarly get 
\begin{align}\label{est-0 critical}
    \dashint_{B_{4r}} \abs{\partial \phi}^d \, \dx \lesssim_{(d)} \dashint_{B_{4r}} f \phi \, \dx.
\end{align}
To estimate the R.H.S. of \eqref{est-0 critical}, we split 
\begin{align*}
    \left| \dashint_{B_{4r}} f \phi \, \dx \right| \le \dashint_{B_{4r}} \abs{f} \abs{\phi - (\phi)_{B_{4r}}} \, \dx + \dashint_{B_{4r}} \abs{f} \abs{(\phi)_{B_{4r}}} \, \dx.
\end{align*}
From \eqref{A and tilde A}, recall the $N$-functions $A(s) := (1+s)\log(1+s) -s$ and $\tilde{A}(s) := \exp(s) -s -1$ where $s>0$. Applying the generalized H\"{o}lder's inequality involving the Orlicz space (see \cite[8.11, Page 269]{AF2003}) with the conjugate pair $(A(s), \tilde{A}(s))$ we estimate the first integral of the above inequality as
\begin{align}\label{gholder}
  \dashint_{B_{4r}} \abs{f} \abs{\phi - (\phi)_{B_{4r}}} \, \dx \le  \frac{1}{|B_{4r}|} \norm{f}_{L^{ A }(B_{4r})} \norm{\phi - (\phi)_{B_{4r}}}_{L^{\tilde{A}}(B_{4r})}. 
\end{align}
Now, we estimate $\norm{\phi - (\phi)_{B_{4r}}}_{L^{\tilde{A}}(B_{4r})}$. For that, we set
\begin{align}\label{M}
    M := \omega_{d} \left( \int_{B_{4r}} \abs{\partial  \phi}^d \right)^{\frac{1}{d}}.
\end{align}
Using the H\"{o}lder's inequality with the conjugate pair $(d, d')$ we get 
\begin{align*}
    \int_{B_{4r} \cap B_{\rho}(x_0)} \abs{\partial  \phi} \le M \rho^{d-1}, \text{ for every } B_{\rho}(x_0) \subset \rd.
\end{align*}
In view of $\tilde{A}$, for any $s>0$, we have 
\begin{align*}
    \int_{B_{4r} } \tilde{ A } \left( \frac{\abs{  \phi - (\phi)_{B_{4r}}  }}{s} \right) \le \int_{B_{4r} } \exp \left(  \frac{\abs{ \phi - (\phi)_{B_{4r}} }}{s} \right).
\end{align*}
For $M$ as in \eqref{M} and for $\sigma_0 = \sigma_0(d)$, we choose $\sigma = \frac{\sigma_0}{2} | B_{4r} | (\diam B_{4r})^{-d}$ and use Proposition \ref{John-Nirenberg} to get 
\begin{align*}
    \int_{ B_{4r} } \exp \left( \frac{\sigma}{M} \abs{ \phi - (\phi)_{B_{4r}} } \right) \lesssim_{(d)} ( \diam B_{4r} )^d \le 1,
\end{align*}
where the last inequality follows from the choice of $r$. Therefore, for $s= \frac{M}{\sigma}$, using the definition of the Orlicz space (Definition \ref{Orlicz}), we conclude that 
\begin{align*}
    \norm{\phi - (\phi)_{B_{4r}}}_{L^{\tilde{A}}(B_{4r})} \le \frac{ \omega_{d} }{ \sigma } \left( \int_{B_{4r}} \abs{\partial  \phi}^d \right)^{\frac{1}{d}} = \frac{ \omega_{d} }{ \sigma } |B_{4r}|^{\frac{1}{d}} \left( \dashint_{B_{4r}} \abs{\partial  \phi}^d \right)^{\frac{1}{d}}.
\end{align*}
Hence from \eqref{gholder} we have 
\begin{align*}
    \dashint_{B_{4r}} \abs{f} \abs{\phi - (\phi)_{B_{4r}}} \le \frac{2 \omega_d}{ \sigma_0 } \abs{B_{4r}}^{\frac{1}{d} - 1} (\text{diam}(B_{4r}))^d \norm{f}_{L^{ A }(B_{4r})} \left( \dashint_{B_{4r}} \abs{\partial  \phi}^d \right)^{\frac{1}{d}},
\end{align*}
where $\abs{B_{4r}}^{\frac{1}{d} - 1} (\text{diam}(B_{4r}))^d = C(d) r^{(d-1)d + 1-d} \lesssim_{(d)} r^{1-d}$, since $r \le 1$. Applying Young's inequality, 
\begin{equation}\label{est critic 1}
    \begin{split}
        \dashint_{B_{4r}} \abs{f} \abs{\phi - (\phi)_{B_{4r}}} & \lesssim_{(d)} r^{1-d}  \norm{f}_{L^{ A }(B_{4r})} \left( \dashint_{B_{4r}} \abs{\partial  \phi}^d \right)^{\frac{1}{d}}  \\
        & \le \ga \dashint_{B_{4r}} \abs{ \partial  \phi }^d \, \dx + C( d, \ga) \left( r^{ 1 - d } \norm{f}_{L^{ A }(B_{4r})} \right)^{\frac{d}{d-1}}.
    \end{split}
\end{equation}
Further, using the H\"{o}lder's and Poincar\'{e} inequalities we estimate 
\begin{equation}\label{est critic 2}
    \begin{split}
        \dashint_{B_{4r}} \abs{f} \abs{(\phi)_{B_{4r}}} = \left( \dashint_{B_{4r}} r \abs{f}  \right) \left( \dashint_{B_{4r}} \frac{\abs{ \phi }}{r} \right) 
        & \le \left( \dashint_{B_{4r}} r \abs{f}  \right) \left( \dashint_{B_{4r}} \left( \frac{\abs{\phi}}{r} \right)^{d}  \right)^{\frac{1}{d}}  \\
        & \lesssim_{(d)} \left( \dashint_{B_{4r}} r \abs{f}  \right)  \left( \dashint_{B_{4r}} 
     \abs{ \partial  \phi }^d  \right)^{\frac{1}{d}}  \\
       & \le \tilde{ \gamma } \dashint_{B_{4r}} \abs{ \partial  \phi }^d  + C( \tilde{\gamma}, d ) \left( \dashint_{B_{4r}} r \abs{f}  \right)^{\frac{d}{d-1}}.
    \end{split}
\end{equation}
Taking $\ga, \tilde{\ga} < \frac{1}{4}$ we obtain from \eqref{est-0 critical}, \eqref{est critic 1}, and \eqref{est critic 2} that 
\begin{align*}
    \dashint_{B_{4r}} \abs{ \partial  \phi }^d & \lesssim_{(d)} \left( \dashint_{B_{4r}} r \abs{f}  \right)^{\frac{d}{d-1}} +  \left( r^{ 1 - d } \norm{f}_{L^{ A }(B_{4r})} \right)^{\frac{d}{d-1}}  \\
    & \lesssim_{(d)} \left( \dashint_{B_{4r}(x_0)}  r \abs{ f } \log \left( e + \frac{\abs{f}}{\int_{B_r(x_0)} \abs{f} \, \dx} \right) \, \dx  \right)^{\frac{d}{ d-1 }},
\end{align*}
where the last inequality holds using Proposition \ref{equivalence} and the fact that $\log(\cdot)$ is an increasing function. Now for a function $h \in W^{1,p}(B_{r}(x_0))$ satisfying the following equation weakly: 
\begin{equation*}
    \begin{aligned}
        - \De_p h &= 0 \mbox{ in } B_{r}(x_0),\\
        h&=v \mbox{ on }  \pa B_{r}(x_0),
    \end{aligned}
\end{equation*}
we proceed with the same arguments as given in Step 2 of Proposition \ref{comparison estimates d>p} to get 
\begin{align*}
    \dashint_{B_r}  \abs{\partial  (v-h)}^p \, \dx \le T r^{p \beta},
\end{align*}
where $\beta = \beta(d,s, p) \in (0,1)$, and $T = T(d,s, \norm{\partial  u}_{L^d(\Omega)}, \norm{u}_{W^{s,d}(\rd)})$. Thus the required estimate \eqref{h est-3} holds. 
\end{proof}

Finally, we prove the comparison estimate in the case $d<p$.

\begin{proposition}\label{comparison estimates for d<p}
Let $d < p$. Suppose $u$ is a weak solution to \eqref{main-eqn}. Then there exists $h$ satisfying $- \De_p h = 0 \mbox{ in } B_r(x_0)$, such that the following estimate holds:
\begin{align}\label{h est-4}
     \dashint_{B_r(x_0)} \abs{ \partial (u - h)(x) }^p \, \dx \le C \left( \dashint_{B_{4r}(x_0)} \abs{r f(x)} \, \dx  \right)^{p'} + T r^{p \be},
 \end{align}
where $C= C(d, p) > 0$ and $T = T(d,s,p, \norm{\partial  u}_{L^p(\Omega)}, \norm{u}_{W^{s,p}(\rd)}) > 0$ are constants, and $\be = \be(d,s,p) \in (0,1)$. 
\end{proposition}

\begin{proof}
    For $v$ as given in \eqref{v-eqn} using Proposition \ref{Morrey} we estimate 
  \begin{align*}
     \left| \dashint_{B_{4r}} f(x) (u-v)(x) \, \dx \right| & \le \sup_{x  \in B_{4r}(x_0)} \abs{(u-v)(x)} \dashint_{B_{4r}} \abs{ f(x) } \, \dx \\
     & = \sup_{x  \in B_{4r}(x_0)} \abs{(u-v)(x)} r^{-1} \dashint_{B_{4r}} \abs{ r f(x) } \, \dx \\
     & \le C(d, p) \left(  \dashint_{ B_{4r} } \abs{ \partial  (u-v)(x) }^p \, \dx \right)^{\frac{1}{p}} \left( \dashint_{B_{4r}} \abs{ r f(x) } \, \dx \right).
  \end{align*}
Further, Young's inequality yields
\begin{align*}
    C(d, p) \left(  \dashint_{ B_{4r} } \abs{ \partial  (u-v) }^p \right)^{\frac{1}{p}} \left( \dashint_{B_{4r}} \abs{ r f(x) }  \right) \le \ga \dashint_{B_{4r}} \abs{ \partial  (u-v) }^p \, \dx +C( d, p , \ga) \left( \dashint_{B_{4r}} \abs{r f(x)} \right)^{p'}.
\end{align*}
Now by taking $\ga < \frac{1}{2}$ rest of the proof holds following the same arguments as given in Proposition \ref{comparison estimates d>p}.
\end{proof}

By considering Proposition \ref{comparison estimates d>p}, Proposition \ref{comparison estimates for d=p}, and Proposition \ref{comparison estimates for d<p}, we have the following result in this section: 

\begin{theorem}\label{comparison estimate final}
    Let $u$ be a weak solution to \eqref{main-eqn}. When $d=p$, we choose $r>0$ small enough so that  $C(d) ( \diam B_{4r} )^d \le 1$, for $C(d)$ as given in Proposition \ref{John-Nirenberg}.  If $(f, \la)$ satisfies the finiteness condition \eqref{H}, then there exists a function $h$ satisfying \[- \De_p h = 0 \mbox{ in } B_r(x_0),\] such that the following estimate holds:
    \begin{align*}
        \left( \dashint_{B_r(x_0)} \abs{ \partial (u - h)(x) }^p \, \dx \right)^{\frac{1}{p}} \le C F(f, B_{4r}(x_0)) \left( \left( \frac{\la}{K_1} \right)^{\max\{ 0, 2-p \}} + \left( \frac{\la}{K_2} \right)^{ \max\{ 0, 2-p \} }\right) + T r^{\be},
    \end{align*}
    where $C = C(d, p) > 0$ and $T = T(d,s,p, \norm{\partial  u}_{L^p(\Omega)}, \norm{u}_{W^{s,p}(\rd)}) > 0$ and $\be = \be(d,s,p) \in (0,1)$.
\end{theorem}

\section{Proof of Theorem \ref{thm:main}}\label{section-4}

In this section, we prove the borderline Lipschitz regularity following the strategy developed in \cite{KM2014MA, KM2014JEMS} by proving certain decay estimates for the solution. As a consequence, by an application of the Lebesgue differentiation theorem, we obtain bounds for all Lebesgue points of the gradient of the solution. We begin by fixing a Lebesgue point $X \in \Omega$ of $\partial u$. Let $0<\rho<1$. We work in the reference ball $B_{10\rho}(X) \subset \Omega$. Additionally, for $d=p$, we choose $\rho >0$ small enough so that  $C(d) ( \diam B_{4 \rho} )^d \le 1$, for $C(d)$ as given in Proposition \ref{John-Nirenberg}. Next, we consider $\sigma = \sigma(d,p) \in (0, \frac{1}{2})$ as the decay parameter to be fixed during the proof (see \eqref{sigma choice}). We now choose universal constants $k, H_1, H_2, H_3 \geq 1$ as follows. With $L = L(d,p) \geq 1$ and $\al = \al(d,p) \in (0,1)$ as in \eqref{eq:8.16} 
we first choose 
\begin{equation}\tag{$\kappa_1$}
\label{k1}
    k \geq 1+\frac{1}{\al}\left(10+\frac{d}{p}\right)
\end{equation}
large enough so that the following holds: 
\begin{equation}\tag{$\kappa_2$}
\label{k2}
     4 L\left(\frac{\sigma}{2}\right)^{-\frac{d}{p}+\al k} \leq \frac{1}{256}.
\end{equation}
For $k$ given above and for $C_0=C_0(d,p) \geq 1$ as in \eqref{eq:8.15} we choose
\begin{equation}\tag{$H_i$}\label{H123}
    \begin{split}
        & H_1 \geq (3)(64)\left(\frac{\sigma}{2}\right)^{-\frac{d}{p}}, \\
        & H_2 \geq \left(2^{100}C_0\left(\frac{\sigma}{2}\right)^{-\frac{d}{p}(k+1)}\right)^{\max\{1,\frac{1}{p-1}\}}, \text{ and } \\
        & H_3 \geq 2^{100}C_0\left(\frac{\sigma}{2}\right)^{-\frac{d}{p}(k+1)},
    \end{split}
\end{equation}
to be some universal constants. Finally, we define $\lambda$ as follows
\begin{equation}\tag{\textbf{$\tilde{\Phi}$}}\label{eq:lambda}
    \begin{split}
        \la^p := H_1^p \left( \dashint_{B_{\rho}(X)} \left( \abs{\partial  u(x)}^p  + 1 \right) \, \dx \right) + H_2^{p} M(f, B_{4\rho}(X))^{\max\{\frac{p}{p-1},p\}} \\ + H_3^p T^p r^{ p \beta } \left(\frac{1}{1 - \left( \frac{\sigma}{2} \right)^{\beta}} \right)^p,
    \end{split}
\end{equation}
where $T > 0$ and $\be = \be(d,s,p)\in (0, 1)$ are given in Theorem \ref{comparison estimate final}. Note that $(f,\la)$ satisfies the finiteness condition \eqref{H} with $K_1 = H_1$ and $K_2 = H_2$, because $(f,\la)$ satisfies \eqref{eq:lambda}.

We now fix a sequence of balls which converge to $X$. For $j \geq 1$ let 
\begin{align*}
    B_j = B_{\rho_j}(X), \text{ where } \rho_j = \left(\frac{\sigma}{2}\right)^{j-1}\rho.
\end{align*}
Note that $B_1$ here stands for the ball of radius $\rho$ centred at $X$. Recalling Definition \ref{defn:pot} we denote
\begin{align*}
    F_j := F(f, B_{4\rho_j}(x_0)).
\end{align*}
We now use a switching radius argument to force the required estimate. To this end, we let
\begin{align*}
    G_j := \left(\dashint_{B_j}|\partial u|^p\,dx\right)^{\frac{1}{p}} + \left(\frac{\sigma}{2}\right)^{-\frac{d}{p}} \mathfrak{E}_p(\partial u, B_j).
\end{align*}
By noticing that $\mathfrak{E}_p(\partial u, B_1) \le 2 (\dashint_{B_1}|\partial u|^p\,dx )^{\frac{1}{p}}$, and $\sigma < 1$, we get 
\begin{align*}
    G_1 \leq 3\left(\frac{\sigma}{2}\right)^{-\frac{d}{p}}\left(\dashint_{B_1}|\partial u|^p\,dx\right)^{\frac{1}{p}}.
\end{align*}
The choice of $\lambda$ in \eqref{eq:lambda} and $H_1$ in \eqref{H123} imply that
\begin{align*}
    G_1 \leq 3\left(\frac{\sigma}{2}\right)^{-\frac{d}{p}}8^{\frac{d}{p}}\frac{\la}{H_1} \le \frac{1}{64}\la.
\end{align*}
Recalling \eqref{k2} we may assume that there exists a step $j_0 \geq 1$ such that 
\begin{align}\label{eqG}
    G_{j_0} \leq \frac{1}{64}\lambda \; \text{ and } \; G_j > \frac{1}{64}\lambda \text{ for } j>j_0,
\end{align}
because if $\eqref{eqG}$ does not hold, then along some subsequence $\underset{j \ra \infty}{\lim} G_j \leq \frac{\la}{64}$ which under the choice of $\lambda$ in \eqref{eq:lambda} already implies Theorem \ref{thm:main}. So without loss of generality, we assume that $\eqref{eqG}$ holds. Next, we establish that even under the assumption \eqref{eqG}, the choices made in \eqref{k1}, \eqref{H123}, and \eqref{eq:lambda} still force the required conclusion. The proof follows from an inductive argument applied to the following claim:
    \begin{paragraph}{\textbf{Claim:}} Suppose there exists $j_1 \geq j_0$ such that for all $j_0 \leq j \leq j_1$ the following inequality is satisfied:
    \begin{equation}\label{eq:indhp}
        \left(\dashint_{B_j}|\partial u|^p\,dx\right)^{\frac{1}{p}} + \mathfrak{E}_p(\partial u, B_j) \leq \lambda.
    \end{equation}
    Then at $j_1+1$'st step, the following estimate holds 
    \[
    \left(\dashint_{B_{j_1+1}}|\partial u|^p\,dx\right)^{\frac{1}{p}} + \mathfrak{E}_p(\partial u, B_{j_1+1}) \leq \lambda.
    \]
    \end{paragraph}
    \begin{paragraph}{Proof of Claim:} The proof is divided into two steps. In the linear settings, the first step is not required since the excess decay condition is direct for the harmonic functions. 
    
    \noi \textbf{Step 1:} In this step, we want to show that the hypothesis \eqref{eq:non-deg} holds for the $p-$harmonic replacement $h_j$ for $u$ in $B_j$, so that the excess decay condition \eqref{eq:exdec} holds for $h_j$. Let $j_0 \leq j \leq j_1$ and $h_j$ be the $p$-harmonic replacement for $u$ in $B_j$. Then the comparison estimates in Theorem \ref{comparison estimate final} imply that 
        \begin{equation}\label{eq:8.15}
            \begin{split}
                & \left(\dashint_{B_j}|\partial u - \partial h_j|^p \,dx\right)^{\frac{1}{p}} \le C_0\left(\left(\frac{\lambda}{H_1}\right)^{\max\{0,2-p\}} + \left(\frac{\lambda}{H_2}\right)^{\max\{0,2-p\}}\right)F_j + C_0T\rho_j^{\beta} \\ 
                & \le C_0 \left(\left(\frac{\lambda}{H_1}\right)^{\max\{0,2-p\}} + \left(\frac{\lambda}{H_2}\right)^{\max\{0,2-p\}}\right)\sum_{j=1}^{\infty}F_j + C_0T\sum_{j=1}^{\infty}\rho_j^{\beta}\\
                & \le C_0 \left(\left(\frac{\lambda}{H_1}\right)^{\max\{0,2-p\}} + \left(\frac{\lambda}{H_2}\right)^{\max\{0,2-p\}}\right)\left(\frac{\lambda}{H_2}\right)^{\min\{p-1,1\}} + C_0T\rho^{\beta}\frac{1}{1-(\frac{\sigma}{2})^{\beta}}\\
                & \le C_0 \left[ \left(\left(\frac{\lambda}{H_1}\right)^{\max\{0,2-p\}} + \left(\frac{\lambda}{H_2}\right)^{\max\{0,2-p\}}\right)\left(\frac{\lambda}{H_2}\right)^{\min\{p-1,1\}} +  \frac{\lambda}{H_3} \right],
            \end{split}
        \end{equation}
        where $C_0 = C_0(d, p) \ge 1$. In the last two lines of \eqref{eq:8.15}, we use Lemma \ref{Dini sum expression}, the condition \eqref{eq:lambda} and the definition of $\rho_j$. 
        Since $H_1,H_2,H_3 \geq 1$, the induction hypothesis \eqref{eq:indhp} and $\eqref{eq:8.15}$ yield
        \begin{align*}
            \left(\dashint_{B_j}|\partial h_j|^p\right)^{\frac{1}{p}} \leq \left(\dashint_{B_j}|\partial u - \partial h_j|^p\right)^{\frac{1}{p}} + \left(\dashint_{B_j}|\partial u|^p\right)^{\frac{1}{p}} \lesssim_{(d,s,p)} \lambda.
        \end{align*}
        It now follows from Proposition \ref{prop:decay1} that
        \begin{align}\label{eq:8.16}
            \sup_{\frac{1}{2}B_j} \,|\partial h_j| \leq L\lambda \; \text{ and } \; \underset{\frac{t}{2}B_j}{\text{osc}} \, \partial h_j \leq Lt^{\alpha} \lambda,
        \end{align}
        for some $L = L(d,p) \ge 1$ and $\alpha = \alpha(d,p) \in (0,1)$ and any $t \in (0,1)$. Thus the upper bound in \eqref{eq:non-deg} is verified with $ v = h_j$. Next, we proceed to verify the lower bound in \eqref{eq:non-deg}. Let $k \geq 1$ be as in \eqref{k1} and \eqref{k2}. Notice that 
        \begin{align*}
            \mathfrak{E}_p(\partial u, B_{j+k}) \le \dashint_{B_{j+k}} \abs{\partial u - (\partial h_j)_{B_{j+k}}}^p \, \dx +  \dashint_{B_{j+k}} \abs{ (\partial h_j)_{B_{j+k}} - (\partial u)_{B_{j+k}} }^p \, \dx,
        \end{align*}
        where using Jensen's inequality, we have
         \begin{align*}
            \abs{ (\partial h_j)_{B_{j+k}} - (\partial u)_{B_{j+k}} }^p = \left| \dashint_{B_{j+k}} \abs{ \partial h_j} - \dashint_{B_{j+k}} \abs{ \partial u}
            \right|^p & \le \left| \dashint_{B_{j+k}} \abs{ \partial u - \partial h_j} \, \dx \right|^p \\
            & \le \dashint_{B_{j+k}} \abs{ \partial u - \partial h_j}^p \, \dx.
         \end{align*} 
      From the above estimates, we compute
        \begin{equation}\label{eq:8.17}
            \begin{split}
            \frac{1}{64}\lambda
            & \overset{\eqref{eqG}}{\leq} 
            G_{j+k} \leq  
            \left(\dashint_{B_{j+k}}|\partial h_j|^p\right)^{\frac{1}{p}} + \left(\frac{\sigma}{2}\right)^{-\frac{d}{p}} \left(\dashint_{B_{j+k}}|\partial u - \partial h_j|^p\right)^{\frac{1}{p}} \\
            & \quad + \left(\frac{\sigma}{2}\right)^{-\frac{d}{p}}\left(\dashint_{B_{j+k}}|\partial u - (\partial h_j)_{B_{j+k}}|^p\right)^{\frac{1}{p}} \\
            & \leq 
            \left(\dashint_{B_{j+k}}|\partial h_j|^p\right)^{\frac{1}{p}} + 3\left(\frac{\sigma}{2}\right)^{-\frac{d}{p}}\left(\dashint_{B_{j+k}}|\partial u - \partial h_j|^p\right)^{\frac{1}{p}} \\
            & \quad + 2\left(\frac{\sigma}{2}\right)^{-\frac{d}{p}}\left(\dashint_{B_{j+k}}|\partial h_j - (\partial h_j)_{B_{j+k}}|^p\right)^{\frac{1}{p}},
            \end{split}
        \end{equation}
      where in the last inequality of \eqref{eq:8.17} comes applying the triangle inequality. From \eqref{eq:8.15} we have  
        \begin{equation*}
        \begin{split}
            & 3\left(\frac{\sigma}{2}\right)^{-\frac{d}{p}}\left(\dashint_{B_{j+k}}|\partial u - \partial h_j|^p\right)^{\frac{1}{p}} \\
            & \leq
            3\left(\frac{\sigma}{2}\right)^{-\frac{d}{p}(k+1)}C_0 \left[ \left(\left(\frac{\lambda}{H_1}\right)^{\max\{0,2-p\}} + \left(\frac{\lambda}{H_2}\right)^{\max\{0,2-p\}}\right)\left(\frac{\lambda}{H_2}\right)^{\min\{p-1,1\}} +  \frac{\lambda}{H_3} \right]. 
        \end{split}
            \end{equation*}
        Also for $t = 2 \left( \frac{ \sigma }{2} \right)^{k}$,  $\eqref{eq:8.16}$ yields 
        \begin{equation*}
          2\left(\frac{\sigma}{2}\right)^{-\frac{d}{p}}\left(\dashint_{B_{j+k}}|\partial h_j - (\partial h_j)_{B_{j+k}}|^p\right)^{\frac{1}{p}} \leq 4 L\left(\frac{\sigma}{2}\right)^{-\frac{d}{p}+\al k} \lambda.
        \end{equation*}
        Further, the choice of $k\geq 1$ in \eqref{k2} implies
        \[
        4 L\left(\frac{\sigma}{2}\right)^{-\frac{d}{p}+\al k} \lambda \leq \frac{1}{256}\lambda.
        \] 
        Moreover, the choices of $H_1, H_2, H_3 \ge 1$ in \eqref{H123}, and $k\geq 1$ in \eqref{k1} imply
       \begin{equation*}
           \begin{split}
               3\left(\frac{\sigma}{2}\right)^{-\frac{d}{p}(k+1)}C_0 \left[ \left(\left(\frac{\lambda}{H_1}\right)^{\max\{0,2-p\}} + \left(\frac{\lambda}{H_2}\right)^{\max\{0,2-p\}}\right)\left(\frac{\lambda}{H_2}\right)^{\min\{p-1,1\}} +  \frac{\lambda}{H_3} \right]  \leq \frac{1}{256}\lambda.
           \end{split}
       \end{equation*}
        In particular, we get
        \[
        \frac{1}{128}\lambda \leq \left(\dashint_{B_{j+k}}|\partial h_j|^p\right)^{\frac{1}{p}} \leq \sup_{B_{j+k}}|\partial h_j| \leq \sup_{B_{j+1}}|\partial h_j|. 
        \]
        Thus the lower bound in \eqref{eq:non-deg} is also verified for $h_j$. Therefore, applying Proposition \ref{prop:decay2}, for any $\gamma \in (0,1)$ there exists $\sigma = \sigma(\gamma,d,p) \in (0, \frac{1}{2})$ so that the following excess decay holds:
        \begin{equation}\label{eq:8.20}
            \mathfrak{E}_p(\partial h_j,B_{j+1}) \leq \gamma \mathfrak{E}_p(\partial h_j, B_j).
        \end{equation}

        \noi \textbf{Step 2:} In this step, using the excess decay condition in \eqref{eq:8.20} we specify the decay $\gamma$, which fixes the decay parameter $\sigma$ universally. Then we prove the claim using the prescribed decay with the large choices of the parameters $H_1, H_2$ and $H_3$ that are made in \eqref{H123}. From the definition of $\mathfrak{E}_p$, and applying Jensen's inequality we have
        \begin{align*}
            & \mathfrak{E}_p(\partial u,B_{j+1}) \leq  2\left(\dashint_{B_{j+1}}|\partial u - \partial h_j|^p\right)^{\frac{1}{p}} + 2\mathfrak{E}_p(\partial h_j,B_{j+1}), \text{ and } \\
            & \mathfrak{E}_p(\partial h_j, B_{j}) \leq 2\left(\dashint_{B_{j}}|\partial u - \partial h_j|^p\right)^{\frac{1}{p}} + 2\mathfrak{E}_p(\partial u, B_j).
        \end{align*}
    Now we use the above estimates and $\eqref{eq:8.20}$ to get 
    \begin{equation*}
            \begin{split}
                \mathfrak{E}_p(\partial u,B_{j+1}) \le 2\left(\frac{\sigma}{2}\right)^{-\frac{d}{p}}\left(\dashint_{B_{j}}|\partial u - \partial h_j|^p\right)^{\frac{1}{p}} + 2 \ga \left( 2\left(\dashint_{B_{j}}|\partial u - \partial h_j|^p\right)^{\frac{1}{p}} + 2\mathfrak{E}_p(\partial u, B_j) \right).
            \end{split}
        \end{equation*}
        Recalling the first inequality of $\eqref{eq:8.15}$ we hence get
        \begin{align*}
            \mathfrak{E}_p(\partial u, B_{j+1}) \le 2C_0\left(\left(\frac{\sigma}{2}\right)^{-\frac{d}{p}} + 2\gamma  \right)\left(\left(\frac{\lambda}{H_1}\right)^{\max\{0,2-p\}} + \left(\frac{\lambda}{H_2}\right)^{\max\{0,2-p\}}\right)F_j \\
            + 2C_0T\left(\left(\frac{\sigma}{2}\right)^{-\frac{d}{p}} + 2\gamma  \right)\rho_j^{\beta}
            + 4\gamma \mathfrak{E}_p(\partial u,B_j).
        \end{align*}
We now take
\begin{align}\label{sigma choice}
    \gamma = \frac{1}{8} \text{ which fixes } \sigma = \sigma (d,p) \in \left( 0, \frac{1}{2} \right), 
\end{align}
 according to Proposition \ref{prop:decay2}. Hence we get
        \begin{align*}
            \mathfrak{E}_p(\partial u, B_{j+1}) \le 2C_0\left(\left(\frac{\sigma}{2}\right)^{-\frac{d}{p}} + 2\gamma  \right)\left(\left(\frac{\lambda}{H_1}\right)^{\max\{0,2-p\}} + \left(\frac{\lambda}{H_2}\right)^{\max\{0,2-p\}}\right)F_j \\
            + 2C_0T\left(\left(\frac{\sigma}{2}\right)^{-\frac{d}{p}} + 2\gamma  \right)\rho_j^{\beta}
            + \frac{1}{2}\mathfrak{E}_p(\partial u,B_j).
        \end{align*}
        Therefore,
        \begin{align*}
            \sum_{j=j_0}^{j_1} \mathfrak{E}_p(\partial u, B_{j+1}) 
            \leq 
          2C_0\left(\left(\frac{\sigma}{2}\right)^{-\frac{d}{p}} + 2\gamma  \right)\left(\left(\frac{\lambda}{H_1}\right)^{\max\{0,2-p\}} + \left(\frac{\lambda}{H_2}\right)^{\max\{0,2-p\}}\right)\sum_{j=j_0}^{j_1}F_j \\
            + 2C_0T\left(\left(\frac{\sigma}{2}\right)^{-\frac{d}{p}} + 2\gamma  \right)\sum_{j=j_0}^{j_1}\rho_j^{\beta} + \mathfrak{E}_p(\partial u, B_{j_0})\\
            \le  2C_0\left(\left(\frac{\sigma}{2}\right)^{-\frac{d}{p}} + 2\gamma  \right)\left(\left(\frac{\lambda}{H_1}\right)^{\max\{0,2-p\}} + \left(\frac{\lambda}{H_2}\right)^{\max\{0,2-p\}}\right)\left(\frac{\lambda}{H_2}\right)^{\min\{p-1,1\}}\\
            + 2C_0\left(\left(\frac{\sigma}{2}\right)^{-\frac{d}{p}} + 2\gamma  \right)\frac{\lambda}{H_3} + \frac{1}{64}\left(\frac{\sigma}{2}\right)^{\frac{d}{p}}\lambda,
        \end{align*}
        where in the last inequality, we use Lemma \ref{Dini sum expression}, the condition \eqref{eq:lambda}, and the assumption on $G_j$ in \eqref{eqG}. In particular, the large  choices of the parameters $H_1, H_2$ and $H_3$ in \eqref{H123} imply
        \begin{equation}\label{eq:8.23}
            \sum_{j=j_0}^{j_1} \mathfrak{E}_p(\partial u, B_{j+1}) \leq \frac{1}{32}\left(\frac{\sigma}{2}\right)^{\frac{d}{p}}\lambda.
        \end{equation}
         Next, we have
        \[
        |(\partial u)_{B_{j+1}}-(\partial u)_{B_{j}}| \leq \dashint_{B_{j+1}}|\partial u - (\partial u)_{B_{j}}| \leq \mathfrak{E}_p(\partial u, B_j) \leq \left(\frac{\sigma}{2}\right)^{-\frac{d}{p}}\mathfrak{E}_p(\partial u, B_j), 
        \]
        which in light of $\eqref{eq:8.23}$ and the assumption on $G_{j_0}$ in \eqref{eqG} implies
        \begin{equation}\label{eq:8.24}
        |(\partial u)_{B_{j_1+1}}| \leq |(\partial u)_{B_{j_0}}| + \left(\frac{\sigma}{2}\right)^{-\frac{d}{p}}\sum_{j=j_0}^{j_1}\mathfrak{E}_p(\partial u, B_j) \leq \frac{1}{16}\lambda.
        \end{equation}
        Finally, using \eqref{eq:8.23}, \eqref{eq:8.24} and $\sigma \in (0, \frac{1}{2})$ we get
        \[
        \left(\dashint_{B_{j_1+1}}|\partial u|^p\,dx\right)^{\frac{1}{p}} + \mathfrak{E}_p(\partial u, B_{j_1+1}) \leq |(\partial u)_{B_{j_1+1}}| + 2\mathfrak{E}_p(\partial u, B_{j_1+1}) \leq \lambda.
        \]
        Thus the claim follows. 
        \end{paragraph}

\section{Acknowledgement}
The authors thank Karthik Adimurthi (TIFR-CAM) for his valuable suggestions, which improved the article. Both authors are grateful for the support provided by the Department of Atomic Energy, Government of India, through project no. 12-R $\&$ D-TFR-5.01-0520.

\bibliographystyle{abbrvnat}

\end{document}